\documentclass[a4paper]{article}
\usepackage{amsmath,amsthm,amsfonts,amssymb}
\usepackage{multicol}
\usepackage{multirow}
\usepackage{diagbox}
\usepackage{graphicx,subfigure}

\usepackage{lineno,cases}
\allowdisplaybreaks 
\usepackage{epstopdf} 
\usepackage{pgf,tikz}
\usetikzlibrary{arrows,shapes,automata,positioning}
\usetikzlibrary{fit}
\usepackage{pst-node}
\usepackage{mathptmx}

\usepackage[pdfstartview=FitH,
CJKbookmarks=true,
bookmarksnumbered=true,
bookmarksopen=true,
colorlinks,
pdfborder=001,
linkcolor=blue,
anchorcolor=blue,
citecolor=blue
]{hyperref}

\usepackage{color,xcolor}
\renewcommand{\top}{{\mathrm{T}}}
\usepackage{bm} 
\newcommand{\lm}{\lambda}
\newcommand{\e}{\varepsilon}
\newcommand{\R}{\mathbb{R}}

\newtheorem{thm}{Theorem}

\newtheorem{lem}{Lemma}
\newtheorem{ass}{Assumption}

\newtheorem{rem}{Remark}
\newtheorem{prop}{Problem}

\newcommand{\pb}{\noindent\textbf{Proof. } }
\newcommand{\pe}{\hfill\rule{4pt}{8pt}}

\def\rm{\mathrm}
\renewcommand{\top}{{\mathrm{T}}}

\begin{document}
	
	\title{Distributed Optimization for a Class of High-order Nonlinear Multi-agent Systems with Unknown Dynamics}
	
	\author{Yutao Tang
		\footnote{This work was supported by National Natural Science Foundation of China under Grant 61503033. Y. Tang is with the School of Automation, Beijing University of Posts and Telecommunications, Beijing 100876, China (e-mail: yttang@bupt.edu.cn). }
	}
	
	\date{}
	
	\maketitle
	
	{\noindent\bf Abstract}: In this paper, we study a distributed optimization problem for a class of high-order multi-agent systems with unknown dynamics. In comparison with existing results for integrators or linear agents, we need to overcome the difficulties brought by the unknown nonlinearities and also the optimization requirement. For this purpose, we employ an embedded control based design and first convert this problem into an output stabilization problem. Then, two kinds of adaptive controllers are given for these agents to drive their outputs to the global optimal solution under some mild conditions. Finally, we show that the estimated parameter vector converges to the true parameter vector under some well-known persistence of excitation condition. The efficacy of these algorithms was verified by a simulation example.


	\section{Introduction}
	
	Over the past few years, distributed optimization of multi-agent systems has become a hot topic due to the fast development of multi-robot networks, machine learning and big data technologies \cite{derenick2007convex,nedic2010constrained, boyd2011distributed}.  In a typical setting of this problem, each agent is assigned with a local cost function and the control objective is to propose distributed controls that guarantee a consensus on the optimal solution of the sum of all local cost functions.  Many effective algorithms have been proposed to achieve this goal in different situations \cite{shi2013reaching,zhu2012distributed,yuan2015gradient,kia2015distributed,yang2016PI,yang2017distributed}. 
	
	Here, we follow this technical line but consider high-order continuous-time nonlinear agents unknown dynamics.  While most of the existing works were only devoted to single-integrator agents, there are many distributed optimization tasks implemented by or depending on physical plants of continuous dynamics in practice, e.g. source seeking in multi-robot systems \cite{zhang2011extremum}, attitude formation control of rigid bodies\cite{song2017relative}, optimal power dispatch over power networks \cite{stegink2017unifying}. Note that these physical dynamics can hardly be described well by single integrators. In fact, an example was given to show that the algorithms designed for single integrators might fail to achieve a desired performance for these high-order agents \cite{tang2018distributed}. Thus, we have to take the high-order and possible nonlinear agents' dynamics into account to achieve the distributed optimization goal.   
	
	In light of the optimization requirement for agents, the gradient-based closed-loop systems are basically nonlinear. Furthermore, the high-order feature of these agents brings many new difficulties to the associated analysis and design, which makes this problem much more challenging.  In fact, very few optimization results have been obtained on this topic with continuous-time agents in the form of high-order dynamics.  For example,  a distributed algorithm for double integrators was proposed in \cite{zhang2017distributed} with an integral control idea and then extended to Euler-Lagrange agents. A class of nonlinear minimum phase agents in output feedback form with unity relative degree was considered by an internal-model design \cite{wang2016cyber}. Distributed optimization problem with bounded controls was also explored for both single and double integrators \cite{xie2017global} . However, distributed optimization of general high-order multi-agent systems is still far from being solved.  
	
	Recently, we proposed an embedded control scheme to solve this problem for general linear systems having well-defined vector relative degrees \cite{tang2017cyb}. To overcome the difficulties brought by both high-order dynamics and possible nonlinearities from gradients, we divided the design into two main parts: optimal signal generator construction and reference tracker design. This embedded technique makes the control design carried out in a ``separative'' way, so as to simplify the whole design by almost independently tackling the optimal consensus problem for single integrators and output tracking problem of high-order agents.
	
	Note that exact information of system matrices is required in this embedded design, which may be not available or with measurement errors in applications. This motivates us to investigate the distributed optimization problem for high-order nonlinear agents with unknown time-varying dynamics. Furthermore, we focus on the case when this unknown time-varying nonlinearity can be linearly parameterized.  The objective of this paper is to extend existing embedded control results to this class of uncertain high-order nonlinear agents and achieve the distributed optimization goal.
	
	In view of the aforementioned observations, the contribution of this paper is at least two-fold. Firstly, a distributed optimization problem was formulated and solved for a group of high-order nonlinear agents, which can be taken as an extended version of existing results for single integrators \cite{nedic2010constrained, kia2015distributed}.  Secondly, the embedded control technique proposed in \cite{tang2017cyb} was further explored and extended to solve the distributed optimization problem of high-order multi-agent systems with unknown dynamics. By removing the requirement of knowing agents' exact dynamics, this work can be taken as an adaptive extension to existing results for linear multi-agent systems \cite{kia2015distributed, xie2017global,tang2017cyb} and nonlinear ones with unity relative degree \cite{wang2016cyber}. Moreover, the obtained conclusions can be applied to achieve an output average consensus for these uncertain nonlinear agents, while only integrators or linear agents were considered in existing works \cite{ren2008distributed,rezaee2015average,tang2017kyb}. 
	
	The remainder of this paper is organized as follows. Problem formulation is presented in Section \ref{sec:formulation}. Then the main result is presented in Sections \ref{sec:main} along with both stability analysis and parameter convergence. Following that, an example is given to illustrate the effectiveness of our algorithms in Section \ref{sec:simu}. Finally, concluding remarks are given in Section \ref{sec:con}.
	
	\textsl{Notations:} Let $\R^n$ be the $n$-dimensional Euclidean space. For a vector $x$, $||x||$ denotes its Euclidean norm. ${\bf 1}_N$ (and ${\bf 0}_N$) denotes an $N$-dimensional all-one (and all-zero) column vector. $\mbox{diag}\{b_1,\,{\dots},\,b_n\}$ denotes an $n\times n$ diagonal matrix with diagonal elements $b_i,\,(i=1,\,{\dots},\,n)$. $\mbox{col}(a_1,\,{\dots},\,a_n) = [a_1^\top,\,{\dots},\,a_n^\top]^\top$ for column vectors $a_i\; (i=1,\,{\dots},\,n)$.  A weighted directed graph (or digraph) $\mathcal {G}=(\mathcal {N}, \mathcal {E}, \mathcal{A})$ is defined as follows, where $\mathcal{N}=\{1,{\dots},n\}$ is the set of nodes, $\mathcal {E}\subset \mathcal{N}\times \mathcal{N}$ is the set of edges, and $\mathcal{A}\in \mathbb{R}^{n\times n}$ is a weighted adjacency matrix \cite{godsil2001algebraic}. $(i,j)\in \mathcal{E}$ denotes an edge leaving from node $i$ and entering node $j$. The weighted adjacency matrix of this digraph $\mathcal {G}$ is described by $A=[a_{ij}]\in \R^{n\times n}$, where $a_{ii}=0$ and $a_{ij}\geq 0$ ($a_{ij}>0$ if and only if there is an edge from agent $j$ to agent $i$).  A path in graph $\mathcal {G}$ is an alternating sequence $i_{1}e_{1}i_{2}e_{2}{\cdots}e_{k-1}i_{k}$ of nodes $i_{l}$ and edges $e_{m}=(i_{m},\,i_{m+1})\in\mathcal {E}$ for $l=1,\,2,\,{\dots},\,k$. If there exists a path from node $i$ to node $j$ then node $i$ is said to be reachable from node $j$. The neighbor set of agent $i$ is defined as $\mathcal{N}_i=\{j\colon (j,\,i)\in \mathcal {E} \}$ for $i=1,\,\dots,\,n$.  A graph is said to be undirected if $a_{ij}=a_{ji}$ ($i,\,j=1,\,{\dots},\,n$). An undirected graph is said to be connected if there is a path between any two vertices.  The weighted Laplacian $L=[l_{ij}]\in \mathbb{R}^{n\times n}$ of graph $\mathcal{G}$ is defined as $l_{ii}=\sum\nolimits_{j\neq i}a_{ij}$ and $l_{ij}=-a_{ij} (j\neq i)$.

	\section{Problem formulation}\label{sec:formulation}
	
	In this paper, we consider a collection of heterogeneous high-order nonlinear systems described by:
	\begin{align}\label{sys:agent}
	\begin{split}
	\dot{x}_{j,i}&=x_{j+1,i}\\ 
	\dot{x}_{n_i,i}&=\Delta_i(x_i,\,\theta_i,\,t)+u_i\\ 
	y_i&=x_{1,i},\quad i=1,\,\dots,\,N,\,j=1,\,\dots,\,n_i-1
	\end{split}
	\end{align}
	where $x_{j,i} \in \R $ is the $j$-th state variable of agent $i$, $x_i\triangleq \mbox{col}(x_{1,i},\,\dots,\,x_{n_i,i}) \in \R^{n_i}$, $y_i \in \mathbb{R}$ and $u_i \in \mathbb{R}$ are respectively the state, output, and input of agent $i$. The function $\Delta_i(x_i,\, \theta_i,\,t)$ represents the unknown time-varying nonlinearities with uncertain parameter $\theta_i=\mbox{col}(\theta_{1,i},\,\dots,\,\theta_{n_{\theta_{n_i}},i})\in \R^{n_{\theta_i}}$, which might result from modeling errors or external perturbations.

	Associated with this multi-agent system, each agent is endowed with a differentiable local cost function $f_i\colon \R \to \R$. The global cost function is defined as the sum of local costs, i.e., $f(y)=\sum\nolimits_{i=1}^{N} f_i(y)$. Moreover, we assume the local cost function $f_i(\cdot)$ is only known to agent $i$ itself and cannot be shared globally in the multi-agent network.  Coupled with these nonlinear agents, we aim to design proper controllers such that the outputs of these agents asymptotically minimize the global cost function.
	
	To clarify our following design, we focus on a class of systems in the form of \eqref{sys:agent} as follows.
	\begin{ass}\label{ass:system}
		For any $i=1,\,\dots,\,N$, there exist a known basis function vector ${\bm p}_i(x_i,\,t)$ and an unknown parameter vector $\theta_i\in \R^{n_{\theta_i}}$ satisfying $\Delta_i(x_i,\,\theta_i,\,t)=\theta_i^\top {\bm p}_{i}(x_{i},\,t)$ for all $x_i\in \R^{n_i}$, $t\geq 0$ and $\theta_i\in \R^{n_{\theta_i}}$.  Furthermore, ${\bm p}_i(x_i,\,t)$ can be uniformly bounded by smooth functions of $x_i$.
	\end{ass}
	
	\begin{rem}
		The unknown nonlinear dynamics of all agents are assumed to be linearly parameterized, which have been widely studied in both classical adaptive control and multi-agent coordination literature \cite{krstic1995nonlinear,yu2012adaptive,hu2014adaptive,wen2017neuro}. 	In fact, a plenty of practical systems can be put into this form, which is general enough to cover integrators, Van der Pol systems, Duffing equations and many mechanical systems.  Moreover, this time-varying feature can further be used to represent many typical external disturbances, e.g.  constants and sinusoidal signals. 
	\end{rem}

	The following assumption is often made in convex optimization literature \cite{bertsekas2003convex,kia2015distributed,qiu2016distributed}, which guarantees the existence and uniqueness of the optimal solution to problem \eqref{opt:main}.
	
	\begin{ass}\label{ass:convexity-strong}
		For $i=1,\,\dots,\,N$, the function $f_i(\cdot)$ is $\underline{l}$-strongly convex and its gradient is $\overline l$-Lipschitz for  constants $\underline{l},\, \overline l>0$.
	\end{ass}
	
	As usual, we assume this optimal solution is finite and denote it as $y^*$, i.e.
	\begin{align}\label{opt:main}
	y^*={\arg}\min_{ y \in \R} \; f(y)=\sum\nolimits_{i=1}^{N} f_i(y)
	\end{align}
	
	Due to the privacy of local cost function $f_i(\cdot)$,  no agent can gather enough information to determine the global optimal solution $y^*$ by itself. Hence, our problem cannot be solved without cooperation and information sharing among these agents. 
	
	For this purpose, we use a weighted undirected graph $\mathcal {G}=(\mathcal {N}, \mathcal {E}, \mathcal{A})$ to describe the information sharing topology with note set $\mathcal{N}=\{1,\,{\dots},\,N\}$. An edge $(i,\,j)\in \mathcal{E}$ between nodes $i$ and $j$ means that agent $i$ and agent $j$ can share information with each other. Suppose the following assumption is fulfilled as in many publications \cite{kia2015distributed,tang2015distributed,zhang2017distributed,tang2017cyb,yang2017distributed}.
	\begin{ass}\label{ass:graph}
		The graph $\mathcal{G}$ is connected.
	\end{ass}
	
	This assumption is about the connectivity of information sharing graph $\mathcal{G}$, which guarantees that any agent's information can reach any other agents. Under this assumption, it is well-known that the associated Laplacian $L$ is positive semidefinite with null space spanned by ${\bm 1}_N$ \cite{godsil2001algebraic}.
	
	The distributed optimization problem considered in this paper is readily formulated as follows.
	\begin{prop}
		For given agents of the form \eqref{sys:agent}, information sharing graph $\mathcal{G}$ and local cost function $f_i(\cdot)$, if possible, determine a distributed protocol for each agent by using its own local data and exchanged information with its neighbors such that
		\begin{itemize}
			\item[1)] all the trajectories of agents are bounded over the time interval $[0,\,+\infty)$.
			\item[2)] their outputs of agents satisfy 
			\begin{align}\label{def:ooc}
			\lim\limits_{t\to +\infty}||y_i(t)-y^*||=0,\quad i=1,\,\dots,\,N
			\end{align}
		\end{itemize}
	\end{prop}
	
	\begin{rem}
		The formulated problem can be taken as a combination of the well-studied topics: distributed optimization \cite{nedic2010constrained, boyd2011distributed, kia2015distributed,yi2015distributed} and output consensus \cite{ren2008distributed, fax2004information, xi2012output}. Since an output consensus of the whole high-order multi-agent system must be achieved as the solution of a convex optimization problem  in the formulated problem, it is certainly more challenging than the existing output consensus results for high-order agents. 
	\end{rem}
	\begin{rem}
		In contrast with existing distributed optimization works, the agents considered here are high-order and heterogeneous, while only single integrators are considered in  many publications \cite{shi2013reaching, yang2016PI,lin2016distributed,qiu2016distributed}. Moreover, we study the case when these agents are with unknown dynamics, while the exact information of agents' dynamics is required in \cite{zhang2017distributed} and \cite{tang2017cyb}. 
	\end{rem}
	
	Particularly, when the local cost function is chosen as $f_i(s)=(s-y_i(0))^2$, this formulation can solve an output average consensus problem for these high-order uncertain nonlinear agents and thus includes exiting results for integrators and linear systems as special cases \cite{ren2008distributed,rezaee2015average,tang2017kyb}.

	As mentioned above, the main difficulty to solve the distributed optimization for these agents lies in the coupling of the high-order structure associated with agents' dynamics and the global optimization requirement. To overcome this point, we adopt the embedded control scheme \cite{tang2017cyb} and propose adaptive control laws to solve the distributed optimization problem for agent \eqref{sys:agent} in the following section.
	
	\section{Main Results}\label{sec:main}
	In this section, we first employ an embedded control approach to convert our problem into an output stabilization problem, and then propose distributed adaptive algorithms for these high-order nonlinear agents with unknown dynamics to achieve the optimization goal along with parameter convergence analysis.
	
	\subsection{Embedded Control Design}
	The embedded control approach was first proposed by Tang et al \cite{tang2017cyb} to solve the distributed optimization problem for high-order linear agents. In this approach, an optimal signal generator must be constructed first by considering the same optimization problem \eqref{opt:main} for a group of single-integrator agents, in order to asymptotically reproduce the optimal solution $y^*$ by a signal $r_i$.  Then, by taking $r_i$ as an output reference signal for agent $i$, this generator is embedded in the feedback loop via some proper interfaces for original agents. In this way, the distributed optimization problem for general linear agents is divided into two simpler subproblems, i.e., construction of an optimal signal generator for single integrators and design of proper trackers for linear agents with output references, which can be independently solved in a modular way.

	The first subproblem is essentially a conventional distributed optimization problem for single integrators and has been well-studied in existing literature \cite{kia2015distributed,yang2016PI,yi2015distributed}. To solve this problem, the following optimal signal generator was proposed for problem \eqref{opt:main} in \cite{tang2017cyb}. 
	\begin{align}\label{sys:generator}
	\begin{split}
	\dot{r}_i&=-\nabla f_i(r_i)-\sum\nolimits_{j=1}^N a_{ij}(\lm_i-\lm_j)\\
	\dot{\lm}_i&=\sum\nolimits_{j=1}^N a_{ij}(r_i-r_j)
	\end{split}
	\end{align}
	
	The effectiveness of \eqref{sys:generator} was already proven in  \cite{tang2017cyb}. We repeat it as follows for completeness.
	\begin{lem}\label{lem:generator}
		Suppose Assumptions \ref{ass:convexity-strong} and \ref{ass:graph} hold. Then, along the trajectory of system \eqref{sys:generator},  $r_i(t)$ approaches the optimal solution $y^*$ exponentially as
		$t\to \infty$ for $i=1,\,\ldots,\,N$.
	\end{lem}
	
	Suppose we know the analytical form of $f_i(\cdot)$ or at least $\nabla f_i(\cdot)$, the optimal signal generator can be implemented independently to generate the minimizer of problem \eqref{opt:main} as showed in this lemma. However, this requirement may cost too much in applications and even impossible in many practical circumstances. Thus, we follow an oracle-based description\cite{nesterov2013introductory} of $f_i(\cdot)$ and consider the case when only the real-time gradient $\nabla f_i(y_i)$ is available upon requests in our following design. 
	
	Furthermore, exact information of system matrices is required for the proposed controllers in \cite{tang2017cyb}. This implies that the corresponding designs can not be directly used to solve our problem for uncertain agents in the form of \eqref{sys:agent}.  Thus, the distributed optimization problem considered in this paper is much more challenging than those addressed for integrators or linear agents in existing works \cite{kia2015distributed,shi2013reaching,qiu2016distributed,yang2016PI,tang2017cyb}.
	
	To deal with these two issues, we adopt a certainty-equivalence design and propose an algorithm as follows:
	\begin{align}\label{ctr:online}
	\begin{split}
	u_i&=-\hat \theta_i^\top {\bm p}_{i}(x_{i},\,t)+\frac{1}{\e^{n_i}}[k_{1i}(x_{1,i}-r_i)+\sum\nolimits_{j=2}^{n_i}\e^{j-1}k_{ji}x_{j,\,i}]\\
	\dot{\hat \theta}_i&=\phi_i(x_i,\,\hat \theta_i,\, r_i,\,t)\\
	\dot{r}_i&=-\nabla f_i(y_i)-\sum\nolimits_{j=1}^N a_{ij}(\lm_i-\lm_j)\\
	\dot{\lm}_i&=\sum\nolimits_{j=1}^N a_{ij}(r_i-r_j) 
	\end{split}
	\end{align} 
	where $\hat \theta_i $ is the estimation of uncertain vector $\theta_i$, and the constants $k_{1i},\,\dots,\,k_{n_i i}$ and function $\phi_i(\cdot)$ to be specified later. Here the constant $\e>0$ is a tunable high-gain parameter to deal with the real-time gradient issue.

	\begin{rem}
		Here, the variables $\hat\theta_i,\,r_i,\,\lm_i$ constitute the local compensator of agent $i$.  It can be verified that this control is indeed distributed in the sense of only using the agents' own local data and exchanged information with their neighbors. Moreover, from its nominal form \eqref{ctr:online}, the above control is composed of two parts, where the last two subsystems can be taken as a modified version of optimal signal generator \eqref{sys:generator} with real-time gradients and the rest compose an adaptive tracking controller such that  $y_i(t)$ can track $r_i(t)$ as $t$ goes to infinity in spite of those unknown dynamics, which confirms the embedded design methodology.
	\end{rem}
	
	Under the above control law, the composite system is then:
	\begin{align}
	\begin{split}
	\dot{x}_{1,\,i}&=x_{2, i} \\ 
	\vdots&\\
	\dot{x}_{n_i,i}&=(\theta_i^\top-\hat \theta_i^\top) {\bm p}_{i}(x_{i},\,t)+\frac{1}{\e^{n_i}}[k_{1i}(x_{1,i}-r_i)+\sum\nolimits_{j=2}^{n_i}\e^{j-1}k_{ji}x_{j,\,i}]\\ 
	\dot{\hat \theta}_i&=\phi_i(x_i,\,\hat \theta_i,\,r_i,\,t)\\
	\dot{r}_i&=-\nabla f_i(y_i)-\sum\nolimits_{j=1}^N a_{ij}(\lm_i-\lm_j)\\
	\dot{\lm}_i&=\sum\nolimits_{j=1}^N a_{ij}(r_i-r_j) 
	\end{split}
	\end{align}
	
	By letting $\hat x_i=\mbox{col}(x_{1,i}-r_i,\,\e x_{2,i},\,\dots,\,\e^{n_i-1}x_{n_i,i})$, it can be further rewritten as:
	\begin{align}\label{sys:closed-spform}
	\begin{split}
	\e \dot{\hat x}_i&=A_i\hat x_i-\e b_{1i} \dot{r}_i+{\e^{n_i}}b_{2i}(\theta_i^\top-\hat \theta_i^\top) {\bm p}_{i}(x_{i},\,t) \\
	\dot{\hat {\theta}}_i&=\phi_i(x_i,\,\hat \theta_i,\,r_i,\,t)\\
	\dot{r}_i&=-\nabla f_i(y_i)-\sum\nolimits_{j=1}^N a_{ij}(\lm_i-\lm_j)\\
	\dot{\lm}_i&=\sum\nolimits_{j=1}^N a_{ij}(r_i-r_j) 
	\end{split}
	\end{align}
	where $b_{1i}=\mbox{col}(1,\,0,\,\dots,\,0)$, $b_{2i}=\mbox{col}(0,\,\dots,\,0,\,1)$ and $A_i=\left[\begin{array}{c|c}
	0&I_{n_{x_i}-1}\\\hline
	k_{1i}&[k_{2i}\,\dots\,k_{n_i\,i}]
	\end{array}\right]$. 
	
	Note that the above system is almost in a singularly perturbed form  except the adaption dynamics $\hat \theta_i$ \cite{khalil2002nonlinear}. By letting $\e=0$, we have $\hat x_i=0$ and $x_{1,i}=r_i$. The resultant quasi-state-state model of the above composite system is exactly the optimal signal generator \eqref{sys:generator}, which in turn guarantees $x_{1,\,i}(t)\to y^*$ as $t$ goes to infinity by Lemma 1.
	
	Based on these observations, our formulated distributed optimization problem for agent \eqref{sys:agent} is converted to a decentralized output stabilization problem for the above system \eqref{sys:closed-spform} with output $\hat x_i$. In other words,  we need to determine proper function $\phi_i(\cdot)$ and constant $\e>0$ such that all trajectories of \eqref{sys:closed-spform} is bounded over $[0,\,+\infty)$ and satisfying $\hat x_{i}(t)\to 0$ as $t$ goes to infinity.
	
	\subsection{Solvability Analysis}
	
	To solve our problem, we denote $\hat x=\mbox{col}(\hat x_1,\,\dots,\,\hat x_N)$,\, $\theta=\mbox{col}(\theta_1,\,\dots,\,\theta_N)$,\, $\hat \theta=\mbox{col}(\hat \theta_1,\,\dots,\,\hat \theta_N)$,\, $r=\mbox{col}(r_1,\,\dots,\,r_N)$,\, $\lm =\mbox{col}(\lm_1,\,\dots,\,\lm_N)$ and $\bar \theta=\theta-\hat \theta$. The whole composite system can be put into a compact form as follows.
	\begin{align}\label{sys:composite-compact}
	\begin{split}
	\dot{\hat x}&=\frac{1}{\e}A\hat  x-B_{1} \dot{r}+EB_{2} {\bm p}^\top (x,\,t) \bar \theta\\
	\dot{\bar {\theta}}&=\phi(x,\,\hat \theta,\,r,\,t) \\
	\dot{r}&=-\nabla \tilde f(y)- L \lm \\
	\dot{\lm}&=L r
	\end{split}
	\end{align}
	where $\tilde f(y)=\sum\nolimits_{i=1}^N f_i(y_i)$,\, $A=\mbox{blockdiag}(A_1,\,\dots,\,A_N)$,\, $B_1=\mbox{blockdiag}(b_{11},\,\dots,\,b_{1N})$,\,$B_2=\mbox{blockdiag}(b_{21},\,\dots,\,b_{2N})$,\,$E=\mbox{blockdiag}(\e^{n_1-1}I_{n_1},\,\dots,\,\e^{n_N-1}I_{n_N})$,\, ${\bm p}(x,\,t)\triangleq \mbox{blockdiag}({\bm p}_1(x_1,\,t),\,\dots,\,{\bm p}_N(x_N,\,t))$,\, $L$ is the Laplacian of $\mathcal{G}$, and the function $\phi(x,\,\hat \theta,\,r,\,t)$ is determined by $\phi_1(\cdot),\,\dots,\,\phi_N(\cdot)$. 
	
	We first choose gain constants $k_{1i},\,\dots,\,k_{n_i\,i}$ such that the polynomial $s^{n_i}-k_{n_i\, i}s^{n_i-1}-k_{2i}s-k_{1i}$ is Hurwitz for any $1\leq i\leq N$. This implies that the following Lyapunov equation 
	\begin{align}
	A_i^\top P_i+P_iA_i=-2I_{n_i}
	\end{align}
	has a unique positive definite solution $P_i$ with compatible dimensions.

	Inspired by existing Lyapunov-based designs \cite{krstic1995nonlinear,hu2014adaptive}, we present the main result of this paper as follows.
	
	\begin{thm}\label{thm:main:online}
		Suppose Assumptions \ref{ass:system}--\ref{ass:graph} hold. Then, the distributed optimization problem determined by \eqref{sys:agent} and \eqref{opt:main} can be solved by controllers of the form \eqref{ctr:online} with $\phi_i(x_i,\,\hat \theta_i,\,r_i,\,t)={\bm p}_{i}(x_{i},\,t) b_{2i}^\top P_i \hat x_i$ for a small enough $\e>0$, where the constant $k_{ji}$ is chosen as above for $i=1,\,\dots,\,N,\,j=1,\,\dots,\,n_i$.
	\end{thm}
	
	\pb  The proof is mainly based on system composition techniques. 
	
	{\it Step 1}: consider the first two subsystems of \eqref{sys:composite-compact}. Let $\hat V_i=\hat W_i+{\e^{n_i-1}}\bar \theta_i^\top \bar \theta_i$ with $\hat W_i=\hat x_i^\top P_i\hat x_i$. Its time derivative along the trajectory of the composite system \eqref{sys:closed-spform} satisfies
	\begin{align*}
	\dot{\hat V}_i= &2\hat x_i^\top P_i [\frac{1}{\e}A_i\hat x_i-b_{1i} \dot{r}_i+{\e^{n_i-1}}b_{2i}\bar \theta_i^\top {\bm p}_{i}(x_{i},\,t)]\\
	&-2\e^{n_i-1}\bar \theta_i^\top\phi_i(x_i,\,\hat \theta_i,\,r_i,\,t)\\
	= & -\frac{2}{\e}\hat x_i^\top\hat x_i-2\hat x_i^\top P_i b_{1i} \dot{r}_i 
	\end{align*}
	
	By Young's inequality, it holds that 
	\begin{align*}
	\dot{\hat V}_i \leq &-\frac{c_1}{\e}||\hat x_i||^2+c_2||\dot{r}_i||^2
	\end{align*}
	for some known positive constants $c_1$ and $c_2$. 
	
	Choose $\hat V=\sum\nolimits_{i=1}^N \hat V_i$. We further have 
	\begin{align}\label{eq:thm1-eq1}
	\dot{\hat V} \leq &-\frac{c_1}{\e}||\hat x||^2+c_2||\dot{r}||^2
	\end{align}

	{\it Step 2}: consider the last two subsystems of \eqref{sys:composite-compact}, which can be rewritten as 
	\begin{align*}
	\dot{r}&=-\nabla \tilde f(r)- L \lm + {\bm h}(r,\,y)\\
	\dot{\lm}&=L r
	\end{align*}
	where the function ${\bm h}(r,\,y)\triangleq \nabla \tilde f(r)-\nabla \tilde f(y)$ is $\overline l$-Lipschitz in $r-y$ by Assumption \ref{ass:convexity-strong}. 
	
	By taking ${\bm h}(r,\,y)$ as perturbations, we let $\mbox{col}(r^*, \lm^*)$ be the equilibrium point of unperturbed $(r,\,\lm)$-system (i.e. when ${\bm h}(r,\,y)\equiv {\bm 0}$).  Note that the unperturbed system is exactly the optimal signal generator, it implies $r^*={\bm 1}_N y^*$ under Assumptions \ref{ass:convexity-strong} and \ref{ass:graph} by Lemma \ref{lem:generator}. In fact, the equilibrium point satisfies $-\nabla f(r^*)-L\lm^*={\bm 0}_N$ and $Lr^*={\bm 0}_{N}$. As a result, there exists a constant $\theta$ such that $r_1^*=\dots=r_N^*=\theta$ since the null space of $L$ is spanned by ${\bm 1}_N$ under Assumption \ref{ass:graph}. Then, one can obtain $\sum\nolimits_{i=1}^N \nabla f_i(\theta)=0$ by ${\bm 1}_N^\top L={\bm 0}_N$, which implies that $\theta$ is an optimal solution of \eqref{opt:main}. By Assumption \ref{ass:convexity-strong}, this implies $\theta=y^*$ and $r^*={\bm 1}_N y^*$.
	
	Letting $\bar r=r-r^*$ and $\bar \lm=\lm-\lm^*$ gives
	\begin{align*}
	\dot{\bar r}&=-{\bm h}(r,\,r^*) -L \bar \lm + {\bm h}(r,\,y)\\
	\dot{\bar \lm}&=L \bar r
	\end{align*}
	
	Inspired by the proof of Lemma \ref{lem:generator} in  \cite{tang2017cyb} , we let $R \in \mathbb{R}^{N \times (N-1)}$ be a matrix satisfying $R^\top {\bm 1}_N={\bm 0}_N$, $R^\top R=I_{N-1}$, and $RR^\top=I_{N}-\frac{1}{N}{\bm 1}_N{\bm 1}_N^\top$. Apparently, the matrix $R$ has a full column rank.  Denote $T=[\frac{{\bm 1}_N}{\sqrt{N}}~R]^\top$ and perform a coordinate transform $\hat \lm=T \bar \lm $. The above system is  equivalent to $\dot{\hat \lm}_1=0$ and
	\begin{align}\label{sys:composite-osg}
	\dot{\bar r}=-{\bm h}(r,\,r^*)-LR \hat \lm_2+{\bm h}(r,\,y),\quad \dot{\hat \lm}_2=R^\top L\bar r
	\end{align}
	Furthermore, the unperturbed system is globally exponential stable at the origin under Assumptions \ref{ass:convexity-strong} and \ref{ass:graph} by Lemma \ref{lem:generator}.
	
	To apply the system composition arguments, one has to investigate the robustness of  system \eqref{sys:composite-osg}.  For this purpose, we let $\hat r\triangleq\mbox{col}(\bar r,\,\hat \lm_2)$ and apply the converse Lyapunov theorem (Theorem 4.15 in \cite{khalil2002nonlinear}) to the unperturbed subsystem, that is, there is a continuously differentiable Lyapunov function $\bar V(\cdot)$ such that
	\begin{gather*}
	c_3||\hat r||^2\leq \bar V(\hat r)\leq c_4 ||\hat r||^2 \\
	\frac{\partial \bar V}{\partial \bar r}[-{\bm h}(r,\,r^*)-LR\hat \lm_2]+\frac{\partial\bar V}{\partial \bar \lm_2}R^\top L\bar r \leq - c_5 ||\hat r||^2 \\
	||\frac{\partial \bar V}{\partial \hat r}||\leq  c_6||\hat r||
	\end{gather*}
	for some positive constants $c_3,\,\dots,\, c_6$.
	
	Along the trajectory of perturbed system \eqref{sys:composite-osg}, one can obtain:
	\begin{align}\label{eq:thm1-eq2}
	\begin{split}
	\dot{\bar V}= &\frac{\partial \bar V}{\partial \bar r}[-{\bm h}(r,\,r^*)-LR\hat \lm_2]+\frac{\partial\bar V}{\partial \bar \lm_2}R^\top L\bar r + \frac{\partial \bar V}{\partial \hat r} {\bm h}(r,\,y)\\
	&\leq - c_5||\hat r||^2+ c_6\bar l||\hat r||||\hat x||
	\end{split}
	\end{align}
	where we use the Lipschitzness of ${\bm h}(r,\,y)$ in $r-y$ and thus in $\hat x$.
	
	{\it Step 3}: consider the stability of the composite system composed by the first two subsystems of \eqref{sys:composite-compact} and system \eqref{sys:composite-osg}. Let $V=\hat V+ c \bar V$ with $c>0$ to be specified later. By using equalities \eqref{eq:thm1-eq1} and \eqref{eq:thm1-eq2},  the derivative of $V$ with respect to $t$ along the whole composite system satisfies
	\begin{align*}
	\dot{V}\leq &-\frac{c_1}{\e}||\hat x||^2+c_2||\dot{r}||^2- c_5||\hat r||^2+ c_6\bar l||\hat r||||\hat x||
	\end{align*}
	Note that $-\nabla \tilde f(r)- L \lm$ and ${\bm h}(r,\,y)$ are both globally Lipschitz in their arguments $\hat r$ and $\hat x$ under Assumption \ref{ass:convexity-strong}. From \eqref{sys:composite-osg}, one can determine a known constant $c_7>0$ satisfying that $||\dot{r}||^2=||\dot{\bar r}||^2  \leq c_7(||\hat r||^2+||\hat x||^2)$ 
	
	By Young's inequality, we have
	\begin{align*}
	\dot{V}\leq &-\frac{c_1}{\e}||\hat x||^2+c_2c_7(||\hat {r}||^2+||\hat x^2||)- cc_5||\hat r||^2 + c^2||\hat x||^2+ c_6^2\bar l^2||\hat r||^2\\
	&\leq -(\frac{c_1}{\e}-c_2c_7-c^2)||\hat x||^2-(cc_5-c_2c_7-c_6^2\bar l^2)||\hat r||^2
	\end{align*}
	Letting $c>\frac{c_2c_7+c_6^2\bar l^2+1}{c_5}$ and $\e<\frac{c_1}{c_2c_7+c^2+1}$ gives that
	\begin{align*}
	\dot{V}\leq & -||\hat x||^2-||\hat{r}||^2
	\end{align*}
	According to the LaSalle-Yoshizawa theorem (Theorem 2.1 in  \cite{krstic1995nonlinear}), we have that all trajectories of the closed-loop system composed of \eqref{sys:composite-compact} and \eqref{sys:composite-osg} are bounded over the time interval $[0,\,+\infty)$ and satisfy that $\lim_{t\to\infty} ||\hat x_i(t)||+||\hat r(t)||=0$. As immediate results, one can conclude the boundedness of $x(t)$, $r(t)$, $\lm(t)$ and $\hat \theta(t)$.  Moreover, we can obtain that $\lim_{t\to\infty} \hat x_{1,i}(t)=0$ and $\lim_{t\to\infty} r_i(t)=y^*$, which implies that
	$$||y_i(t)-y^*||\leq||x_{1,i}-r_i(t)||+||r_i(t)-y^*||\to 0$$ as $t\to +\infty$.  The proof is thus complete. 
	\pe
	
	\begin{rem}
		Note that we consider uncertain high-order agents in the form of \eqref{sys:agent}, which includes integrators as its special cases. Thus, the theorems can be taken as adaptive extensions of existing results with exact known dynamics\cite{zhang2017distributed,rezaee2015average,tang2017kyb,tang2017cyb}. Moreover, many actuating disturbances can be represented in the form of \eqref{sys:agent}, thus we provide different methods to achieve disturbance rejection from the internal model-based approach used in existing works \cite{tang2015distributed,wang2016cyber}.   
	\end{rem}

	Particularly, when the analytical form of $f_i(\cdot)$ or $\nabla f_i(\cdot)$ is known to us, the optimal signal generator can be implemented independently. Following a similar proof, we can choose the gain parameter $\e$ as any positive constant to solve our problem. 
	
	\begin{thm}\label{thm:main:offline}
		Suppose Assumptions \ref{ass:system}--\ref{ass:graph} hold. Then, the distributed optimization problem determined by \eqref{sys:agent} and \eqref{opt:main} can be solved by the following control
		\begin{align}\label{ctr:offline}
		\begin{split}
		u_i&=-\hat \theta_i^\top {\bm p}_{i}(x_{i},\,t)+\frac{1}{\e^{n_i}}[k_{1i}(x_{1,i}-r_i)+\sum\nolimits_{j=2}^{n_i}\e^{j-1}k_{ji}x_{j,\,i}]\\
		\dot{\hat \theta}_i&={\bm p}_{i}(x_{i},\,t) b_{2i}^\top P_i \hat x_i\\
		\dot{r}_i&=-\nabla f_i(r_i)-\sum\nolimits_{j=1}^N a_{ij}(\lm_i-\lm_j)\\
		\dot{\lm}_i&=\sum\nolimits_{j=1}^N a_{ij}(r_i-r_j) 
		\end{split}
		\end{align} 
		where the constant $k_{ji}$ is chosen as above and $\e>0$ is arbitrary for $i=1,\,\dots,\,N,\,j=1,\,\dots,\,n_i$.
	\end{thm}
	
	\begin{rem}
		In some circumstances, we may further let $\phi_i(x_i,\,\hat \theta_i,\,r_i,\,t)=\Lambda_i{\bm p}_{i}(x_{i},\,t) b_{2i}^\top P_i \hat x_i$ with a positive definite matrix $\Lambda_{i}$. This matrix is called an adaption gain in literature\cite{ioannou1995robust}. It can be used to achieve a fast adaption and then improve the transient performance of our controllers to solve the distributed optimization problem.
	\end{rem}
	
	\begin{rem}
		Without further information of the unknown dynamics, the two controllers may fail in practical applications if there are external disturbances or noises in measurements of $x_i$, although it is theoretically proved to achieve the optimization goal as $t$ goes to infinity. To tackle this problem, we can employ a $\sigma$-modification \cite{ioannou1995robust} for $\hat \theta_i$  with sacrificing some accuracy in control performance as follows:
		\begin{align}\label{ctr:sigma}
		\begin{split}
		u_i&=-\hat \theta_i^\top {\bm p}_{i}(x_{i},\,t)+\frac{1}{\e^{n_i}}[k_{1i}(x_{1,i}-r_i)+\sum\nolimits_{j=2}^{n_i}\e^{j-1}k_{ji}x_{j,\,i}]\\
		\dot{\hat \theta}_i&=-\sigma_{\theta_i}\hat\theta_i+{\bm p}_{i}(x_{i},\,t) b_{2i}^\top P_i \hat x_i\\
		\dot{r}_i&=-\nabla f_i(y_i)-\sum\nolimits_{j=1}^N a_{ij}(\lm_i-\lm_j)\\
		\dot{\lm}_i&=\sum\nolimits_{j=1}^N a_{ij}(r_i-r_j) 
		\end{split}
		\end{align} 
		where $\sigma_{\theta_i}>0$ is a tunable parameter such that $\lim_{t\to+\infty}||y_i(t)-y^*||$ can be smaller than any desired positive constant.
	\end{rem}

	\subsection{Parameter Convergence}
	It has been shown that parameter convergence is essential in achieving robustness of the proposed adaptive controllers\cite{ioannou1995robust,krstic1995nonlinear,mazenc2009uniform}.  
	From the proof of Theorem \ref{thm:main:online}, one can merely conclude that the estimator $\hat \theta_i$ converges to some constant. In this subsection, we assert conditions under which $\hat \theta_i(t)$ will converge to its true value $\theta_i$ as $t$ tends to infinity. 
	
	For this purpose, we further assume the basis function ${\bm p}_i(x_i,\,t)$ satisfying the following condition.
	\begin{ass}\label{ass:PE}
		For any $i=1,\,\dots,\,N$, along the trajectory of the closed-loop system composed of \eqref{sys:agent} and \eqref{ctr:online} or \eqref{ctr:offline}, there exist positive constants $m$,\,$t_0$,\,$T_0$ such that the function ${{\bm p}}_i(x_i(t),\,t)$ is uniformly bounded and the following inequality is satisfied:
		\begin{align}\label{eq:PE}
		\frac{1}{T_0}\int_{t}^{t+T_0}{\bm p}_i(x_i(\tau),\,\tau){\bm p}^\top(x_i(\tau),\,\tau){\rm d}\tau \geq mI_{n_{\theta_i}},\quad \forall t\geq t_0
		\end{align}
	\end{ass}
	
	Note that $x_i(t)$ is ultimately bounded by Theorem \ref{thm:main:online}, the boundedness of ${\bm p}_i(x_i(t),\,t)$ is not too strict. The inequality \eqref{eq:PE} is often known as the persistence of excitation (PE) condition, which has widely used in adaptive control literature \cite{krstic1995nonlinear,hu2014adaptive}. 
	
	\begin{thm}\label{thm:parameter}
		Suppose Assumptions \ref{ass:system}--\ref{ass:PE} hold. Then, along the trajectory of system \eqref{sys:agent} under the controllers proposed in Theorems \ref{thm:main:online} and \ref{thm:main:offline} , it holds that $\lim_{t\to+\infty}\hat \theta_i(t)=\theta_i$ for  $i=1,\,\dots,\,N$.
	\end{thm}
	\pb To show this theorem, we first claim that $\lim_{t\to+\infty}\bar \theta_i^\top(t){\bm p}_i(x_i(t),\,t)=0$.  By the proof of Theorem \ref{thm:main:online}, we have $\hat x_i(\infty)=\int_{0}^{+\infty}\dot{\hat x}_i(\tau){\rm d}\tau=0$. From the uniform boundedness of associated variables and Assumption \ref{ass:PE}, it follows that $\ddot{\hat x}_i(t)$ is also bounded. Using Barbalat's lemma (Lemma 8.2 in \cite{khalil2002nonlinear}) to $\dot{\hat x}_i(t)$ implies that $\lim_{t\to+\infty}\dot{\hat x}_i(t)=0$, which confirms this claim.
	
	Next, noting $\dot{\bar \theta}_i=\dot{\hat \theta}_i={\bm p}_{i}(x_{i},\,t) b_{2i}^\top P_i \hat x_i$ gives $\lim_{t\to+\infty}\dot{\bar \theta}_i(t)=0$ by Assumption \ref{ass:system}. According to Lemma 1 in  \cite{ortega1993asymptotic} or its proof, the two facts $\lim_{t\to+\infty}\dot{\bar \theta}_i(t)=0$ and $\lim_{t\to+\infty}\bar \theta_i^\top(t){\bm p}_i(x_i(t),\,t)=0$ provide us that $\lim_{t\to+\infty}\bar \theta_i(t)=0$ under Assumption \ref{ass:PE}. The proof is thus complete. 
	\pe
	
	\begin{rem}
		Since the unknown dynamics is linearly parameterized by Assumption \ref{ass:system}, this theorem can be further modified and applied to any number of components
		in ${\bm p}_i(x_i,\,t)$ satisfying such a PE condition, and then used to address the parameter convergence problem in a more precise way. Specially, when the basis function is time-invariant, the $j$-th component ${\bm p}_{j,i}(x_i)$ of ${\bm p}_i(x_i)$ is persistently excited if $\lim_{x_i\to \mbox{col}(y^*,\,0,\,\dots,\,0)}{\bm p}_{j,i}(x_i)\neq 0$ and guarantees the convergence of $\hat \theta_{j,i}(t)$ to $\theta_{j,i}$ as $t$ goes to infinity.	
	\end{rem}

	\section{Simulations}\label{sec:simu}
	
	In this section, we present a numerical example to illustrate our problem and the effectiveness of our designs.  
	
\begin{figure}
		\centering
		\includegraphics[width=0.80\textwidth]{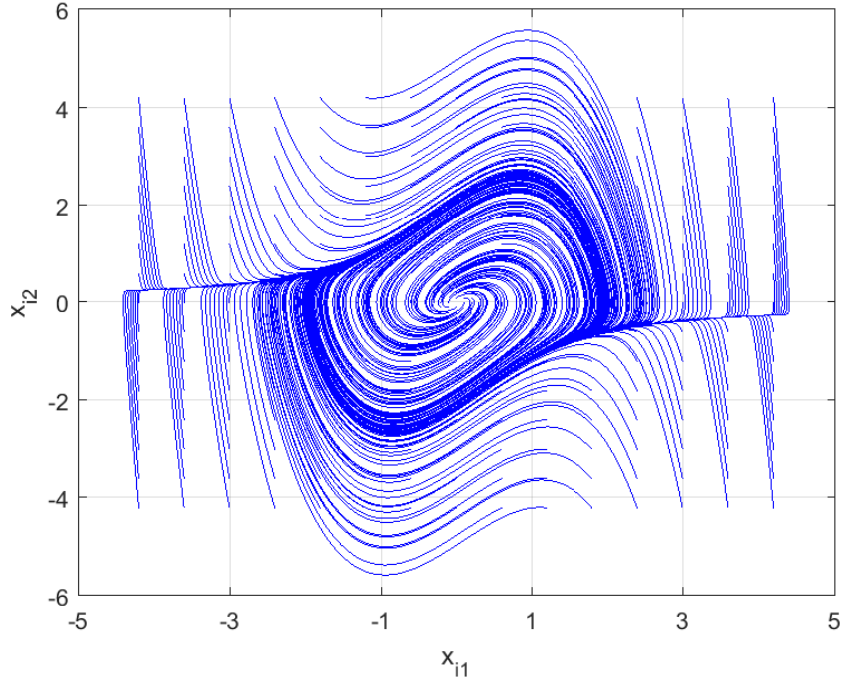}
		\caption{Phase portraits of the unforced Van der Pol system.}\label{fig:phase}
	\end{figure}

	Consider a multi-agent system including four controlled Van der Pol systems as follows.
	\begin{align*}
	\dot{x}_{1,i}&=x_{2, i}\\
	\dot{x}_{2,i}&=\Xi_i+u_i\\
	y_i&=x_{1,i},\quad i=1,\,2,\,3,\,4
	\end{align*}
	where $\Xi_i\triangleq -a_i x_{1,i}+b_i (1-x_{1,i}^2)x_{2,i}$ with $a_i,\,b_i>0$ but unknown. The trajectories of the unforced system with different initial conditions when $a_i=b_i=1$ are depicted in Figure \ref{fig:phase}.

	\begin{figure}
		\centering
		\begin{tikzpicture}[scale=2, shorten >=1pt, node distance=2.2 cm, >=stealth',
		every state/.style ={circle, minimum width=0.8cm, minimum height=0.8cm}]
		\node[align=center,state](node1) {\scriptsize 1};
		\node[align=center,state](node2)[right of=node1]{\scriptsize 2};
		\node[align=center,state](node3)[right of=node2]{\scriptsize 3};
		\node[align=center,state](node4)[right of=node3]{\scriptsize 4};
		\path[-]  (node1) edge (node2)
		(node2) edge (node3)
		(node3) edge (node4)
		(node1) edge [bend left=60]  (node3)
		;
		\end{tikzpicture}
		\caption{Information sharing graph $\mathcal G$ in our example.}\label{fig:graph}
	\end{figure}
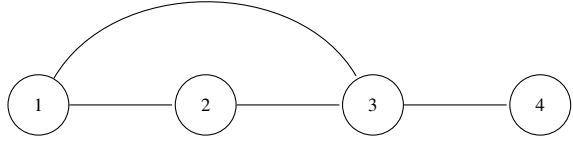

	To make it more interesting, we further assume that agent $i$ is subject to an actuating disturbance $d_i(t)$ described by 
	$$d_i(t)=D_i v_i(t),\quad \dot{v}_i=S_iv_i$$
	where $D_i=[1~0]$ and $S_i=[0~1;\,-1~0]$.
	
	The information sharing graph of this multi-agent system is depicted in Figure \ref{fig:graph} with unity edge weights. The local cost functions are as follows.
	$${f_1}(y) = (y-8)^2,\, {f_2}(y) = \frac{y^2}{{20\sqrt {y^2 + 1} }} + y^2,\, {f_3}(y) =  \frac{y^2}{80\ln {({y^2} + 2} )} + (y -5)^2,\, {f_4}(y) =  \ln \left( {{e^{ - 0.05{y}}} + {e^{0.05{y}}}} \right) + y^2$$
	
	Denote $\Delta_i(x_i,\,t)=\Xi_i+d_i(t)$. The agents are of the form \eqref{sys:agent}. Note that $d_i(t)=A_{1i}\sin(t)+A_{2i}\cos(t)$ for some constants $A_{1i}$,\, $A_{2i}$ depending upon $v_i(0)$. By letting $\theta_i=\mbox{col}(\theta_{1,i},\,\dots,\,\theta_{4,i})=\mbox{col}(a_i,\,b_i,\,A_{1i},\,A_{2i})$ and ${\bm p}_i(x_i,\,t)=\mbox{col}(-x_{1,i},\,(1-x_{1,i}^2)x_{2,i},\,\sin(t),\,\cos(t))$, one can find that Assumption \ref{ass:system} is fulfilled. Moreover, both Assumptions \ref{ass:convexity-strong} and \ref{ass:graph} are also verified. Additionally, the optimal solution is $y^*\approx 3.24$ by numerically minimizing the global cost function $f(y)=\sum\nolimits_{i=1}^4 f_i(y)$. According to Theorem \ref{thm:main:online}, the distributed optimization problem for these agents can be solved by a controller of the form \eqref{ctr:online}.  
	
	\begin{figure}
		\centering
		\includegraphics[width=0.90\textwidth]{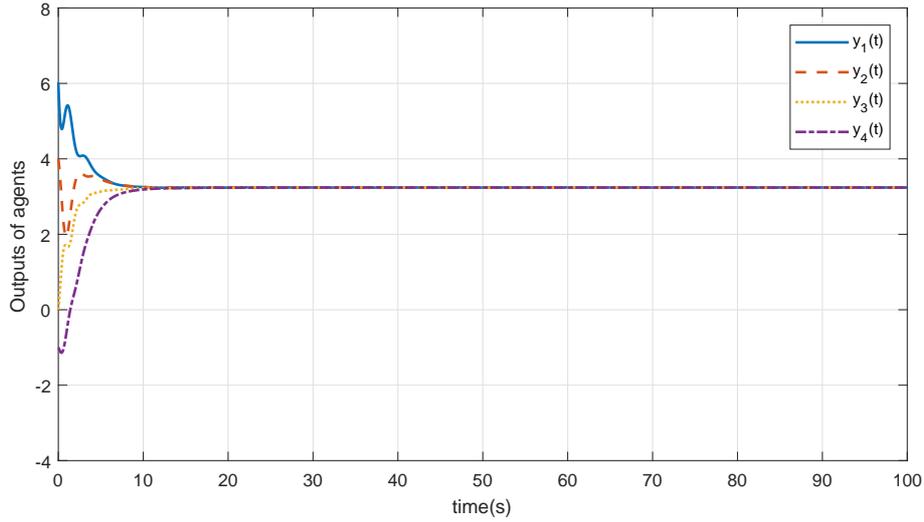}
		\caption{Profiles of $y_i(t)$ under the controller \eqref{ctr:online}.}\label{fig:simu}
	\end{figure}
	
	\begin{figure}
		\centering
		\subfigure[]{
			\includegraphics[width=0.42\textwidth]{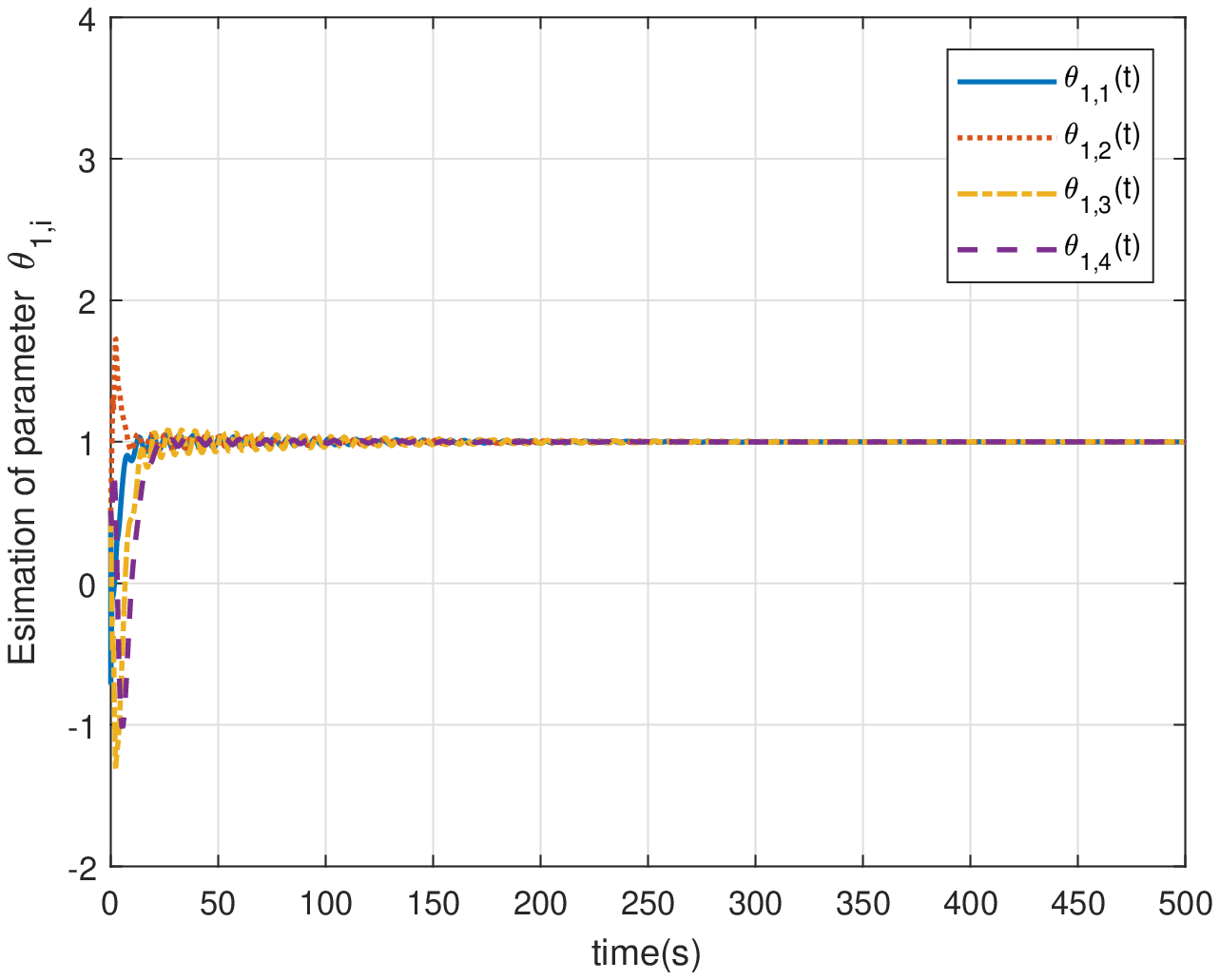}
		}
		\subfigure[]{
			\includegraphics[width=0.42\textwidth]{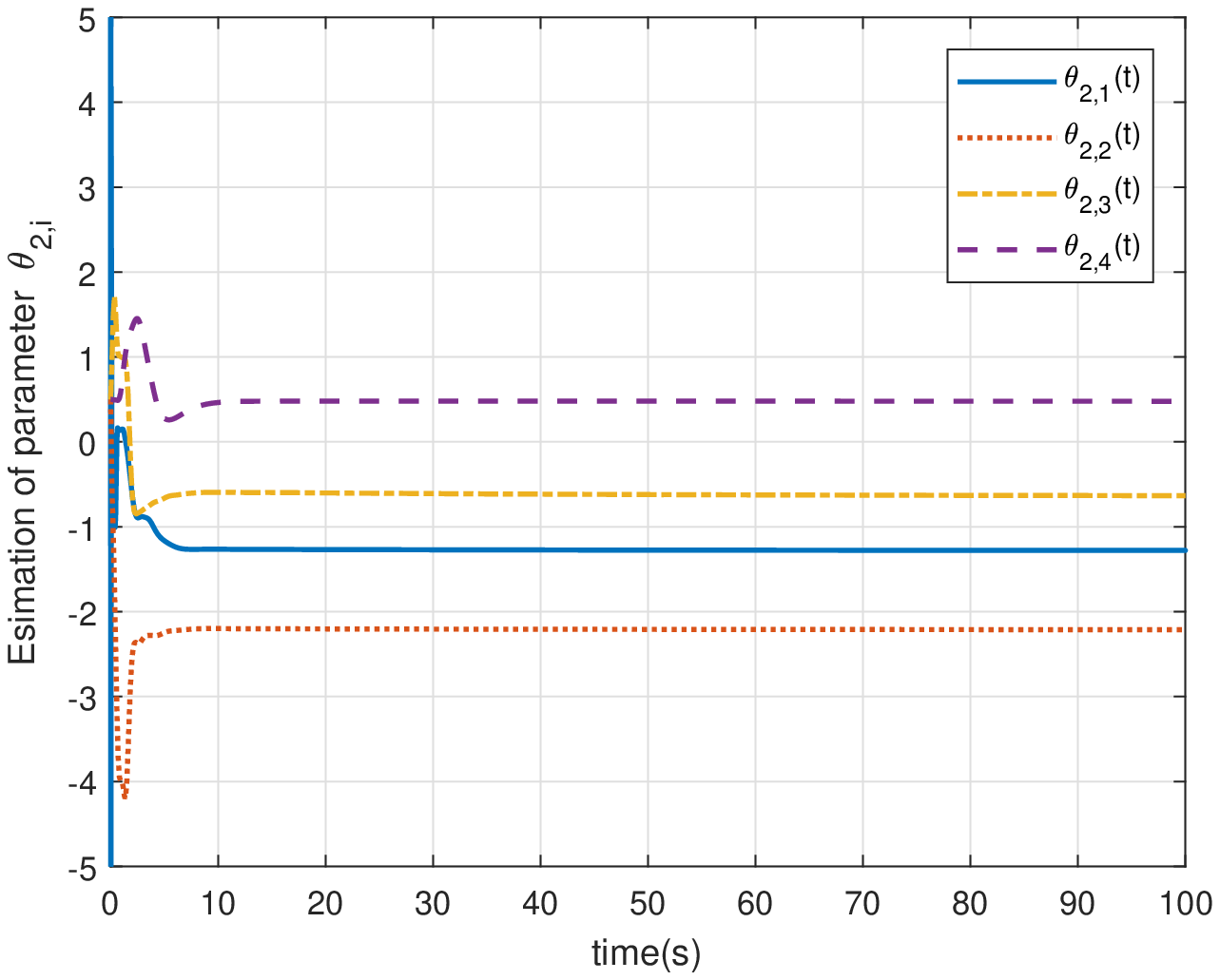}
		}
		\caption{Profiles of $\hat \theta_{1,i}(t)$ and $\hat\theta_{2,i}$ under the controller \eqref{ctr:online}.}\label{fig:parameter-1}
	\end{figure}
	
\begin{figure}
		\centering
		\subfigure[]{
			\includegraphics[width=0.42\textwidth]{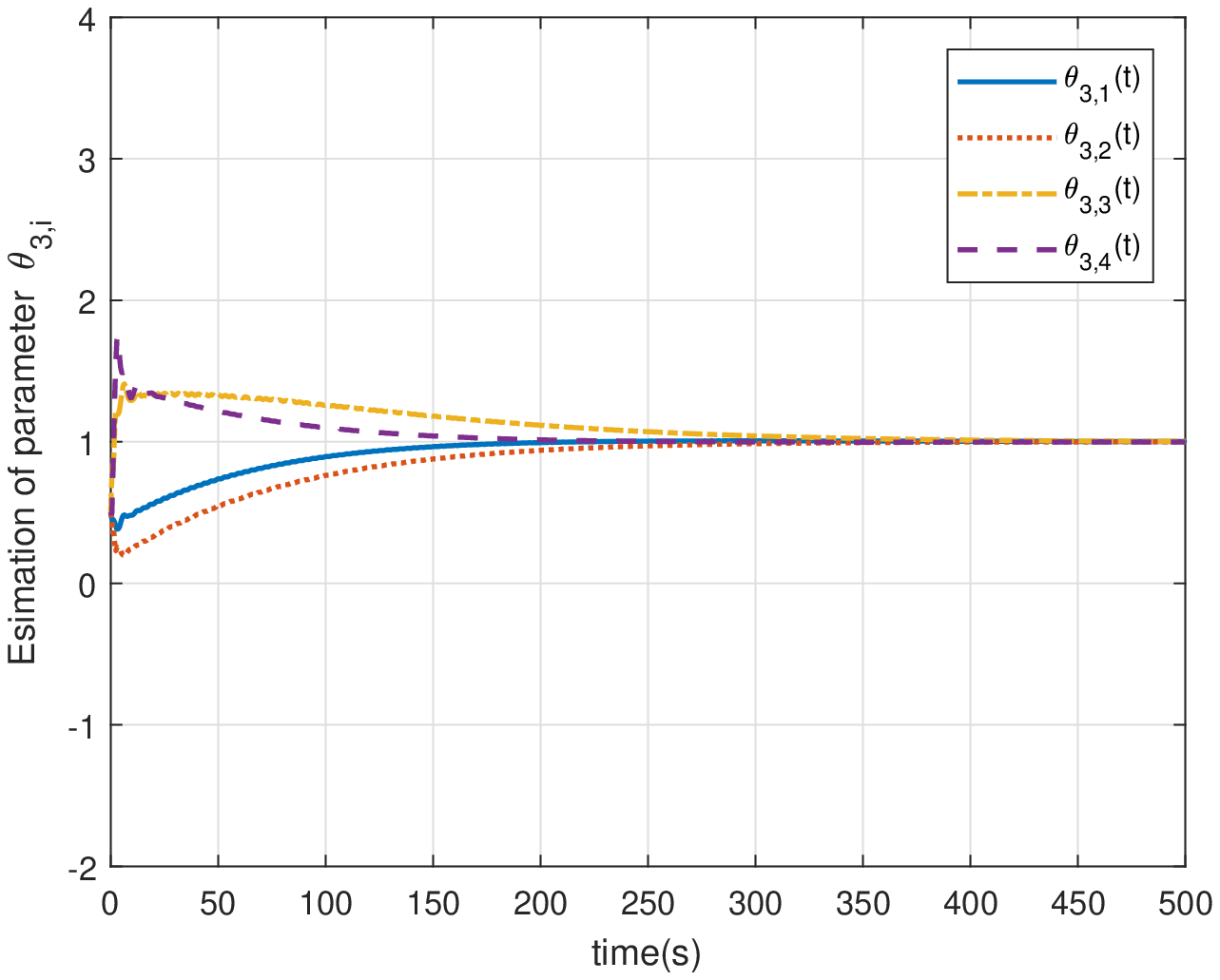}
		}\quad
		\subfigure[]{
			\includegraphics[width=0.42\textwidth]{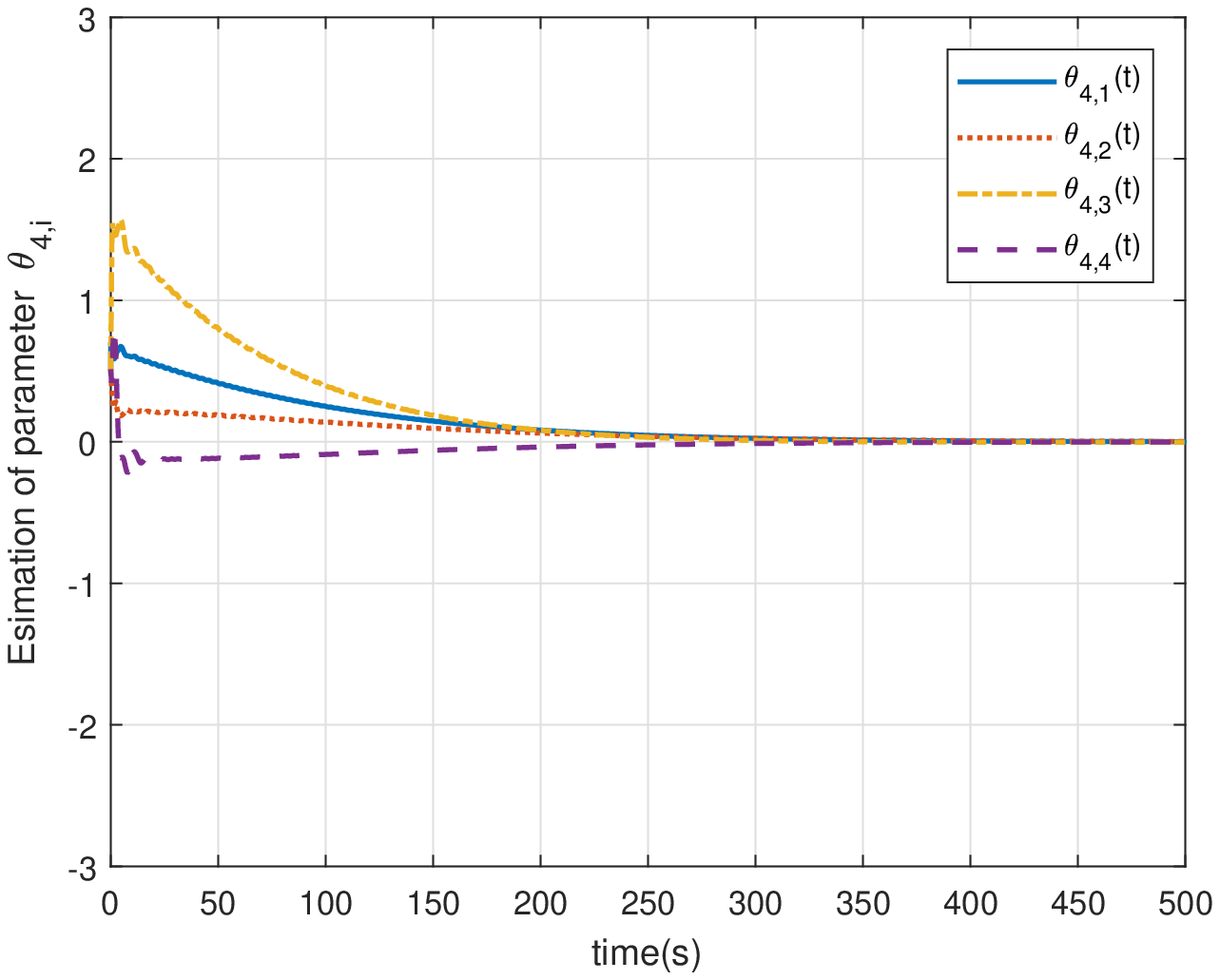}
		}
		\caption{Profiles of $\hat \theta_{3,i}(t)$ and $\hat\theta_{4,i}$ under the controller \eqref{ctr:online}.}\label{fig:parameter-2}
	\end{figure}

	For simulations, we choose the parameters as $a_i=b_i=1$, $v_i(0)=\mbox{col}(1,\,0)$, $k_{1i}=-4$, $k_{2i}=-4$, $\Lambda_{i}=10 I_4$ and $\e=0.2$. The profiles of agents' outputs are shown in Figure \ref{fig:simu}.  Satisfactory performance is observed. For parameter convergence, we have ${\bm p}_{1,i}({\bm x}_i,\,t)=-x_{1,i}$, which is time-invariant and satisfies that $\lim_{x_i\to \mbox{col}(y^*,\,0,\,\dots,\,0)}{\bm p}_{1,i}({\bm x}_i)=-y^*\neq 0$. By some calculations, one can also determine that $$\int_{-\frac{\pi}{2}}^{\frac{\pi}{2}}[{\bm p}_{3,i}(\tau)~~{\bm p}_{4,i}(\tau)]^\top[{\bm p}_{3,i}(\tau)~~{\bm p}_{4,i}(\tau)]{\rm d}\tau=\int_{-\frac{\pi}{2}}^{\frac{\pi}{2}}[\sin^2(\tau)~\sin(\tau)\cos(\tau); \sin(\tau)\cos(\tau)~\cos^2(\tau)]{\rm d}\tau=[\pi/2~~0;\,0~~{\pi/2}]. $$ Using Theorem \ref{thm:parameter}, we can conclude that the estimators $\hat \theta_{1,i}$,\,$\hat \theta_{3,i}$,\, $\hat \theta_{4,i}$ will converge to their true values, while $\hat \theta_{2,i}$ may fail. This conclusion is confirmed by Figures \ref{fig:parameter-1} and \ref{fig:parameter-2}.

	\section{Conclusions}\label{sec:con}
	A distributed optimization problem was formulated for a class of high-order nonlinear systems with unknown dynamics. By using an embedded control scheme, we proposed distributed adaptive controls to solve this problem under standard assumptions and parameter convergence was also addressed. Output feedback control with directed information sharing graphs will be our future work.

	\bibliographystyle{ieeetr}
	\bibliography{opt_adaptive}
	
\end{document}